\documentclass[12pt]{article}
\usepackage[utf8]{inputenc}
\usepackage[margin=1in]{geometry}
\usepackage{amsmath}
\usepackage{amsthm}
\usepackage{amsfonts}
\usepackage{amssymb}
\usepackage{txfonts}
\usepackage[backend=biber, sorting=none, maxnames=4]{biblatex}
\usepackage{ragged2e}
\usepackage{authblk}
\providecommand{\keywords}[1]
{
  \small	
  \textbf{\textit{Keywords: }} #1
}

\providecommand{\subjclass}[1]
{
  \small	
  \textbf{\textit{Mathematics Subject Classification: }} #1
}

\newtheorem{theorem}{Proposition}
\newtheorem{corollary}[theorem]{Corollary}

\title{Certain products of sums of Lambert series}
\author{\emph{Aung Phone Maw}}
\date{}
\begin{document}

\maketitle

\begin{abstract}
     Using elementary means, we prove an identity giving the infinite product form of a sum of Lambert series originally stated by Venkatachaliengar, then rediscovered by Andrews, Lewis, and Liu. Then we derive two identities expressing certain products of sums of Lambert series.
\end{abstract}
\raggedright
. \\
.
\par

\section{Preamble}
Throughout this exposition, we shall use the following standard notation : 
\[ (x;q)_{\infty} = (x)_{\infty} = \prod_{i\geq 1} (1-xq^{i-1}) , \]\[ (x;q)_{n} = (x)_n = \frac{(x;q)_{\infty}}{(xq^n;q)_{\infty}} , \quad n \in \mathbb{C},\quad (x;q)_0 = 1. \]

Then, let us recall, the following two important results from the theory of basic hypergeometric series, the Jacobi triple product identity\cite[p.~18]{nathanJ} and the next is just the partial fraction decomposition of $ \frac{1}{(z)_{\infty} (z^{-1}q)_{\infty}}$ (this identity which appears in \cite[p.~452]{Andrews1984HeckeMF} can be confirmed by studying the behaviour of functions on both sides of the equation at the points $z=q^n,\quad n\in \mathbb{Z}$.) : 
\[ \tag{1.1} \label{eq1} (q)_{\infty} (z)_{\infty} (z^{-1}q)_{\infty} = \sum_{n=-\infty}^{+\infty} (-1)^n z^n q^\frac{n^2 - n}{2} , \]
\[ \tag{1.2} \label{eq2} \frac{{(q)_{\infty}}^2}{(z)_{\infty} (z^{-1}q)_{\infty}}  = \sum_{n=-\infty}^{+\infty} \frac{(-1)^n q^\frac{n^2 + n}{2}}{1-zq^n} .\]

In this exposition, we shall prove three identities on sums of Lambert series. The first being an identity giving the infinite product form of a sum of Lambert series. The equivalent form of this statement was in fact originally stated by Venkatachaliengar\cite[p.~37]{Venkatachaliengar}, in which Venkatachaliengar stated the fundamental multiplicative identity of the Jordan-Kronecker function. It was then rediscovered by Andrews, Lewis and Liu \cite{AndrewsLewisLiu}. For our purpose, this identity will be stated in the following form.

\begin{theorem}

For $a,b,c$ arbitrary, there holds : 
\begin{multline}
    \tag{1.3} \label{eq6} 
    \frac{{(q)_{\infty}}^2 (ab)_{\infty} (a^{-1}b^{-1}q)_{\infty} (ac)_{\infty} (a^{-1}c^{-1}q)_{\infty} (bc)_{\infty} (b^{-1}c^{-1}q)_{\infty}}{(a)_{\infty} (a^{-1}q)_{\infty} (b)_{\infty} (b^{-1}q)_{\infty} (c)_{\infty} (c^{-1}q)_{\infty} (abc)_{\infty} (a^{-1}b^{-1}c^{-1}q)_{\infty}}\\
    = 1+\sum_{r\geq 0} \frac{aq^r}{1-aq^r} - \sum_{r\geq 1} \frac{a^{-1}q^r}{1-a^{-1}q^r} + \sum_{r\geq 0} \frac{bq^r}{1-bq^r} - \sum_{r\geq 1} \frac{b^{-1}q^r}{1-b^{-1}q^r} \\+ \sum_{r\geq 0} \frac{cq^r}{1-cq^r} - \sum_{r\geq 1} \frac{c^{-1}q^r}{1-c^{-1}q^r} + \sum_{r\geq 1} \frac{a^{-1}b^{-1}c^{-1}q^r}{1-a^{-1}b^{-1}c^{-1}q^r} - \sum_{r\geq 0} \frac{abcq^r}{1-abcq^r}.
\end{multline}
\end{theorem}

The following two identities relate certain sums of Lambert series to products of sums of Lambert series. Our main objective will be in deriving them.

\vspace{0.69cm}
\begin{theorem}
    For $b,c$ arbitrary, there holds : 
    \begin{flalign}
        &\notag\left\{ \sum_{r\geq 0} \frac{bq^r}{1-bq^r} - \sum_{r\geq 1} \frac{b^{-1}q^r}{1-b^{-1}q^r}-\sum_{r\geq 0} \frac{bcq^r}{1-bcq^r} + \sum_{r\geq 1} \frac{b^{-1}c^{-1}q^r}{1-b^{-1}c^{-1}q^r} \right\}\left\{ \sum_{r\geq 0} \frac{cq^r}{1-cq^r} - \sum_{r\geq 1} \frac{c^{-1}q^r}{1-c^{-1}q^r} \right. &\\& \notag\left.
        -\sum_{r\geq 0} \frac{bcq^r}{1-bcq^r} + \sum_{r\geq 1} \frac{b^{-1}c^{-1}q^r}{1-b^{-1}c^{-1}q^r} \right\} = \frac{bc}{(1-bc)^2}+\sum_{r \geq 1} r(b^rc^r + b^{-r}c^{-r})\frac{q^r}{1-q^r} &\\& \tag{1.4} \label{eq8}\hspace{7cm}+\sum_{r \geq 1} (b^r+b^{-r}+c^r+c^{-r}-b^rc^r-b^{-r}c^{-r} -2)\frac{q^r}{(1-q^r)^2}.&
    \end{flalign}
\end{theorem}
\vspace{0.69cm}
\begin{theorem}
    For $b,c$ arbitrary, there holds :
    \begin{flalign}
    &\notag\left\{ \frac{1}{2} + \sum_{r\geq 0} \frac{bq^r}{1-bq^r} - \sum_{r\geq 1} \frac{b^{-1}q^r}{1-b^{-1}q^r}+ \sum_{r\geq 0} \frac{cq^r}{1-cq^r} - \sum_{r\geq 1} \frac{c^{-1}q^r}{1-c^{-1}q^r}-\sum_{r\geq 0} \frac{bcq^r}{1-bcq^r} + \sum_{r\geq 1} \frac{b^{-1}c^{-1}q^r}{1-b^{-1}c^{-1}q^r} \right\}^2 &\\&\notag
    =\frac{1}{4} + \sum_{r \geq 0} \left\{ \frac{bq^r}{(1-bq^r)^2} +  \frac{cq^r}{(1-cq^r)^2} +\frac{bcq^r}{(1-bcq^r)^2} \right\} + \sum_{r \geq 1} \left\{ \frac{b^{-1}q^r}{(1-b^{-1}q^r)^2} +  \frac{c^{-1}q^r}{(1-c^{-1}q^r)^2} +\frac{b^{-1}c^{-1}q^r}{(1-b^{-1}c^{-1}q^r)^2} \right\} &\\& \tag{1.5} \label{eq11} \quad-6 \sum_{r \geq 1} \frac{q^r}{(1-q^r)^2} .&
    \end{flalign}
\end{theorem}
\vspace{0.69cm}

 We shall begin our exposition by studying the following function $ f $, which is also known as the Jordan-Kronecker function :
\[ f(a,b;q) = f(a,b) = \sum_{n=-\infty}^{+\infty} \frac{a^n}{1-bq^n} .\]

\section{Some q-Identities on $f(a,b)$}
First, let us write down some trivial facts about $f$ :
\begin{align*}
    \tag{2.1}f(a,b)&=f(b,a), \\
    \tag{2.2}f(a,b)&=-b^{-1} f(qa^{-1},b^{-1}),\\
    \tag{2.3}f(a,b)&=\frac{1-ab}{(1-a)(1-b)} + \sum_{n=1}^{\infty} \frac{ba^n q^n}{1-bq^n} - \sum_{n=1}^{\infty} \frac{b^{-1}a^{-n} q^n}{1-b^{-1}q^n},\\
\end{align*}
\begin{align}
    \tag{2.4} \lim_{a\to1} \left\{f(a,b)-\frac{1}{1-a}\right\} &= \sum_{n=0}^{\infty} \frac{b q^n}{1-bq^n} - \sum_{n=1}^{\infty} \frac{b^{-1} q^n}{1-b^{-1}q^n},\\
    \tag{2.5} a^k f(a,bq^k) &= f(a,b) ,\\
    \tag{2.6} \frac{df(a,b)}{da} = f_a(a,b)&= \sum_{n=-\infty}^{+\infty} \frac{b^n q^n}{(1-aq^n)^2} .
\end{align}
\par
We now consider the following expression : 
\[ (q)_{\infty} (a)_{\infty} (a^{-1}q)_{\infty} f(a,b) ,\]
which we note that, can be resolved in the following manner : 
\begin{flalign}
    \notag&\left\{  \sum_{k=-\infty}^{+\infty} (-1)^k a^k q^\frac{k^2 - k}{2}  \right\} \left\{ \sum_{n=-\infty}^{+\infty} \frac{a^n}{1-bq^n} \right\}&\\
    \notag&= \sum_{k, n = -\infty}^{+\infty} (-1)^k a^{-k} q^\frac{k^2 + k}{2} \frac{a^{n+k}}{1-bq^{n+k}}&\\
    \notag&=\sum_{k, n = -\infty}^{+\infty} (-1)^k q^\frac{k^2 + k}{2} \frac{a^{n}}{1-bq^{n+k}}&\\
    \notag&=\sum_{k, n = -\infty}^{+\infty} (-1)^n a^n b^n q^\frac{n^2 -n}{2} \frac{(-1)^k q^\frac{k^2 + k}{2}}{1-bq^{k}}&\\
    \notag&=\left\{  \sum_{n=-\infty}^{+\infty} (-1)^n a^n b^n q^\frac{n^2 - n}{2}  \right\} \left\{ \sum_{k=-\infty}^{+\infty} \frac{(-1)^k q^\frac{k^2 + k}{2}}{1-bq^k} \right\}&\\
    \tag{2.7} &= \frac{{(q)_{\infty}}^3 (ab)_{\infty} (a^{-1}b^{-1}q)_{\infty}}{(b)_{\infty} (b^{-1}q)_{\infty}} .&
\end{flalign}
Thus, we now see that $f$ has the following product representation :
\[ \tag{2.8} \label{eq3} f(a,b) = \frac{{(q)_{\infty}}^2 (ab)_{\infty} (a^{-1}b^{-1}q)_{\infty}}{(a)_{\infty} (a^{-1}q)_{\infty}(b)_{\infty} (b^{-1}q)_{\infty}} .\]

Then, we shall consider another important expression, which is :
\[ f(a,b)f(a,c) .\]
We see that this expression can also be resolved, in the following way :
\begin{flalign*}
    &\left\{ \sum_{k=-\infty}^{+\infty} \frac{a^k}{1-bq^k} \right\} \left\{ \sum_{n=-\infty}^{+\infty} \frac{a^n}{1-cq^n} \right\}&\\
    &= \sum_{k, n = -\infty}^{+\infty}\frac{a^k a^{n-k}}{(1-bq^k)(1-cq^{n-k})}&\\
    &= \sum_{k, n = -\infty}^{+\infty}\frac{a^{n}bq^k}{(1-bq^k)(bq^k-bcq^{n})}&
\end{flalign*}
\begin{flalign*}
    &= \sum_{k, n = -\infty}^{+\infty}\frac{a^{n}bq^k}{1-bcq^n} \left\{ \frac{1}{bq^k-bcq^n} + \frac{1}{1-bq^k}\right\} &\\
    &= \sum_{n = -\infty}^{+\infty} \frac{a^{n}}{1-bcq^n}\left\{ \sum_{k=0}^{\infty}\left( \frac{1}{1-cq^{n-k}} + \frac{bq^k}{1-bq^k} \right) + \sum_{k=1}^{\infty}\left( \frac{1}{1-cq^{n+k}} + \frac{bq^{-k}}{1-bq^{-k}} \right)  \right\} &\\
    &= \sum_{n = 0}^{\infty} \frac{a^{n}}{1-bcq^n}\left\{ \sum_{k=0}^{\infty}\left( \frac{1}{1-cq^{n-k}} + \frac{bq^k}{1-bq^k} \right) + \sum_{k=1}^{\infty}\left( \frac{cq^{n+k}}{1-cq^{n+k}} + \frac{1}{1-bq^{-k}} \right)  \right\} &\\& \quad + \sum_{n = 1}^{\infty} \frac{a^{-n}}{1-bcq^{-n}}\left\{ \sum_{k=0}^{\infty}\left( \frac{1}{1-cq^{-n-k}} + \frac{bq^k}{1-bq^k} \right) + \sum_{k=1}^{\infty}\left( \frac{cq^{-n+k}}{1-cq^{-n+k}} + \frac{1}{1-bq^{-k}} \right)  \right\} &\\
    &= \sum_{n = 0}^{\infty} \frac{a^{n}}{1-bcq^n}\left\{ \sum_{k=0}^{\infty}\left( \frac{1}{1-cq^{-k}} + \frac{bq^k}{1-bq^k} \right) + \sum_{k=1}^{\infty}\left( \frac{cq^{k}}{1-cq^{k}} + \frac{1}{1-bq^{-k}} \right) + \sum_{k=0}^{n-1} \frac{1}{1-cq^{n-k}} - \sum_{k=1}^{n} \frac{cq^k}{1-cq^k} \right\} &\\& \quad + \sum_{n = 1}^{\infty} \frac{a^{-n}}{1-bcq^{-n}}\left\{ \sum_{k=0}^{\infty}\left( \frac{1}{1-cq^{-k}} + \frac{bq^k}{1-bq^k} \right) + \sum_{k=1}^{\infty}\left( \frac{cq^{k}}{1-cq^{k}} + \frac{1}{1-bq^{-k}} \right) - \sum_{k=0}^{n-1} \frac{1}{1-cq^{-k}} + \sum_{k=1}^{n} \frac{cq^{-n+k}}{1-cq^{-n+k}} \right\} &\\
    &= \sum_{n = 0}^{\infty} \frac{na^n}{1-bcq^n} + \sum_{n = 1}^{\infty} \frac{(-n)a^{-n}}{1-bcq^{-n}} + \left\{ \sum_{n=-\infty}^{+\infty} \frac{a^n}{1-bcq^n} \right\}\left\{ 1+ \sum_{r=0}^{\infty} \frac{bq^r}{1-bq^r} - \sum_{r=1}^{\infty} \frac{b^{-1}q^r}{1-b^{-1}q^r} \right. &\\& \left. \quad + \sum_{r=0}^{\infty} \frac{cq^r}{1-cq^r} - \sum_{r=1}^{\infty} \frac{c^{-1}q^r}{1-c^{-1}q^r} \right\}.
\end{flalign*}
We note that what we just arrived is in fact the fundamental multiplicative identity by Venkatachaliengar\cite[p.~37]{Venkatachaliengar}. And since $\sum_{n = 0}^{\infty} \frac{na^n}{1-bcq^n} + \sum_{n = 1}^{\infty} \frac{(-n)a^{-n}}{1-bcq^{-n}} = af_a(a,bc)$, using the product representation in \eqref{eq3}, we have :
\[ \tag{2.9} \label{eq4} af_a(a,bc) = f(a,bc)\left\{ \sum_{r=0}^{\infty} \frac{aq^r}{1-aq^r} - \sum_{r=1}^{\infty} \frac{a^{-1}q^r}{1-a^{-1}q^r} + \sum_{r=1}^{\infty} \frac{a^{-1}b^{-1}c^{-1}q^r}{1-a^{-1}b^{-1}c^{-1}q^r} - \sum_{r=0}^{\infty} \frac{abcq^r}{1-abcq^r} \right\}.\]
Thus, we arrive at the following result :
\begin{multline}
     \tag{2.10} \label{eq5} f(a,b)f(a,c)=f(a,bc)\left\{ 1+\sum_{r\geq 0} \frac{aq^r}{1-aq^r} - \sum_{r\geq 1} \frac{a^{-1}q^r}{1-a^{-1}q^r} + \sum_{r\geq 0} \frac{bq^r}{1-bq^r} - \sum_{r\geq 1} \frac{b^{-1}q^r}{1-b^{-1}q^r} \right. \\ \left. + \sum_{r\geq 0} \frac{cq^r}{1-cq^r} - \sum_{r\geq 1} \frac{c^{-1}q^r}{1-c^{-1}q^r} + \sum_{r\geq 1} \frac{a^{-1}b^{-1}c^{-1}q^r}{1-a^{-1}b^{-1}c^{-1}q^r} - \sum_{r\geq 0} \frac{abcq^r}{1-abcq^r} \right\} . 
\end{multline}

Finally, using the product representation of $f$, we arrive at \textbf{Proposition 1}. 
Before proceeding any further, let us define the function $l$ as follows :
\[  l(b) = \sum_{r\geq 0} \frac{bq^r}{1-bq^r} - \sum_{r\geq 1} \frac{b^{-1}q^r}{1-b^{-1}q^r} = \lim_{a \to 1} \left\{ f(a,b) - \frac{1}{1-a} \right\} .\]

Next, we shall consider the following expression : 
\[ \left\{ f(a,b) - \frac{1}{1-a} \right\}\left\{ f(a,c) - \frac{1}{1-a} \right\}. \]
Which can be rewritten as :
\begin{flalign*}
    & f(a,b)f(a,c) - \frac{f(a,b)}{1-a} - \frac{f(a,c)}{1-a} + \frac{1}{(1-a)^2}&\\
    &= f(a,bc)\{ 1+l(a)+l(b)+l(c)-l(abc) \} - \frac{f(a,b)}{1-a}- \frac{f(a,c)}{1-a} + \frac{1}{(1-a)^2}&\\
    &= \left\{f(a,bc)-\frac{1}{1-a} + \frac{1}{1-a}\right\}\left\{ 1+l(a)+l(b)+l(c)-l(abc) -\frac{a}{1-a} + \frac{1}{1-a}\right\} &\\& \quad - \frac{f(a,b)}{1-a}- \frac{f(a,c)}{1-a} + \frac{1}{(1-a)^2}&\\
    &= \left\{f(a,bc)-\frac{1}{1-a} \right\}\left\{ l(a)+l(b)+l(c)-l(abc) -\frac{a}{1-a} \right\} + \frac{f(a,bc)-\frac{1}{1-a}}{1-a} &\\& \quad + \frac{l(a)+l(b)+l(c)-l(abc) -\frac{a}{1-a}}{1-a} -\frac{\left( f(a,b)-\frac{1}{1-a} \right)}{1-a}-\frac{\left( f(a,c) - \frac{1}{1-a}\right)}{1-a} &\\
    &= \left\{f(a,bc)-\frac{1}{1-a} \right\}\left\{ l(a)+l(b)+l(c)-l(abc) -\frac{a}{1-a} \right\} &\\
    & \quad + \frac{f(a,bc)-\frac{1}{1-a} - l(abc) + l(b)- \left( f(a,b) - \frac{1}{1-a}\right) + l(c)- \left( f(a,c) - \frac{1}{1-a}\right) + l(a)-\frac{a}{1-a}}{1-a},
\end{flalign*}
now let $a \to 1$, then we get : 
\begin{multline} \tag{2.11} \label{eq7}
    l(b)l(c) = l(bc)\{l(b)+l(c)-l(bc)\} + \sum_{n \geq 1} (b^n+b^{-n}+c^n+c^{-n}-b^nc^n-b^{-n}c^{-n}-2)\frac{q^n}{(1-q^n)^2} \\+ \sum_{n \geq 0} \frac{bcq^n}{(1-bcq^n)^2} + \sum_{n \geq 1} \frac{b^{-1}c^{-1}q^n}{(1-b^{-1}c^{-1}q^n)^2}.
\end{multline}
Thus, we finally arrive at \textbf{Proposition 2} by making the elementary observation which is :  $l(bc)^2-l(bc)\{l(b)+l(c)\}+l(b)l(c) = \{l(bc)-l(b)\}\{l(bc)-l(c)\}$.
Now, if we proceed further and let, $c \to b^{-1}$, we have :
\[ \tag{2.12} \label{eq9}l(b)^2 =  b\frac{dl(b)}{db} - l(b) - 2\sum_{r \geq 1} (b^r+b^{-r}) \frac{q^r}{(1-q^r)^2} - 2 \sum_{r \geq 1} \frac{q^r}{(1-q^r)^2} . \]
Thus, we can conclude from \eqref{eq7} and \eqref{eq9} that : 
\begin{multline}
\tag{2.13} \label{eq10}
\{ \frac{1}{2} + l(b) + l(c) - l(bc) \}^2 = \frac{1}{4} + \sum_{r \geq 0} \left\{ \frac{bq^r}{(1-bq^r)^2} +  \frac{cq^r}{(1-cq^r)^2} +\frac{bcq^r}{(1-bcq^r)^2} \right\} \\+ \sum_{r \geq 1} \left\{ \frac{b^{-1}q^r}{(1-b^{-1}q^r)^2} +  \frac{c^{-1}q^r}{(1-c^{-1}q^r)^2} +\frac{b^{-1}c^{-1}q^r}{(1-b^{-1}c^{-1}q^r)^2} \right\} -6 \sum_{r \geq 1} \frac{q^r}{(1-q^r)^2} .
\end{multline} 
Which then we note is in fact identical to \textbf{Proposition 3}. Now in the next section, we shall unearth some noteworthy particular cases implied by this section.

\section{Particular Cases and Consequences}
In \textbf{Proposition 1}, let $ q \to q^7$, $a=-q$, $b=-q^2$, $c=-q^4$, then we get the following corollary.
\vspace{0.69cm}
\begin{corollary}
\begin{multline}
\tag{3.1} \label{eq12}
    \frac{1}{2} + \sum_{r\geq 0} \left\{ \frac{q^{7r+3}}{1+q^{7r+3}} + \frac{q^{7r+5}}{1+q^{7r+5}} + \frac{q^{7r+6}}{1+q^{7r+6}} -\frac{q^{7r+1}}{1+q^{7r+1}} -\frac{q^{7r+2}}{1+q^{7r+2}} -\frac{q^{7r+4}}{1+q^{7r+4}} \right\}\\= \frac{(q^7;q^7)_\infty (q;q)_\infty}{2(-q^7;q^7)_\infty(-q;q)_\infty}. 
\end{multline}
\end{corollary}

\noindent This identity was originally found by Ramanujan \cite[p.~304]{RamaIII}. Using the identity $\frac{(q;q)_\infty}{(-q;q)_{\infty}} = \sum_{n=-\infty}^{+\infty} (-1)^n q^{n^2}$, restated number theoretically, this is equivalent to the following statement : 
\vspace{0.69cm}
\begin{corollary}
\begin{flalign*} 
\text{\textit{The number of ways a natural number }} N \text{\textit{ can be written in the form }} 7a^2 + b^2\hspace{0.2cm} (a,\hspace{0.2cm} b \in \mathbb{N})\\ = 
\begin{cases}
 \frac{(-1)^N}{2} \left\{ C(N) + (-1)^{N-1} \right\} &\text{\textit{if N or }$\frac{N}{7}$\textit{ is a square,}} \\
        \frac{(-1)^N}{2} C(N)  &\text{\textit{otherwise.}}
\end{cases}\\
\end{flalign*}
\begin{flalign*}
&Where,&\\
&C(N) = \sum_{\substack{d|N \\ d \equiv 1\pmod{7}}} (-1)^\frac{N}{d} - \sum_{\substack{d|N \\ d \equiv 6\pmod{7}}} (-1)^\frac{N}{d} +\sum_{\substack{d|N \\ d \equiv 2\pmod{7}}} (-1)^\frac{N}{d} -\sum_{\substack{d|N \\ d \equiv 5\pmod{7}}} (-1)^\frac{N}{d} &\\& \tag{3.2} \quad \quad \quad + \sum_{\substack{d|N \\ d \equiv 4\pmod{7}}} (-1)^\frac{N}{d} -\sum_{\substack{d|N \\ d \equiv 3\pmod{7}}} (-1)^\frac{N}{d}. &\\
\end{flalign*}
\end{corollary}
\vspace{0.69cm}

Now we continue to look out for other consequences of our main results, in \textbf{Proposition 1}, let $ q \to q^9$, $a=q$, $b=q^2$, $c=q^3$, we get : 
\vspace{0.69cm}
\begin{corollary}
    \begin{flalign}
       \notag 1+\sum_{r\geq 0}\left\{ \frac{q^{9r+1}}{1-q^{9r+1}} - \frac{q^{9r+8}}{1-q^{9r+8}} + \frac{q^{9r+2}}{1-q^{9r+2}}  - \frac{q^{9r+7}}{1-q^{9r+7}} +2\frac{q^{9r+3}}{1-q^{9r+3}} -2\frac{q^{9r+6}}{1-q^{9r+6}} \right\} \\ \tag{3.3}=\frac{(q^9;q^9)_{\infty} ^3(q^4;q^9)_{\infty} ^3(q^5;q^9)_{\infty} ^3}{(q;q)_{\infty}}.
    \end{flalign}
\end{corollary}

Let $q \to q^5, b=c=q,$ in \textbf{Proposition 2}, then we get :
\vspace{0.69cm}
\begin{corollary}
  \begin{multline} \label{eqMod5(a)}
    \left\{ \sum_{r \geq 0} \left( \frac{q^{5r+1}}{1-q^{5r+1}}-\frac{q^{5r+2}}{1-q^{5r+2}}+\frac{q^{5r+3}}{1-q^{5r+3}} -\frac{q^{5r+4}}{1-q^{5r+4}} \right) \right\}^2 \\  = \sum_{r \geq 0} \frac{q^{5r+2}}{(1-q^{5r+2})^2} +  \sum_{r \geq 0} \frac{q^{5r+3}}{(1-q^{5r+3})^2}+ \sum_{r \geq 1} \left( \frac{2rq^{5r+1}}{1-q^{5r+1}} - \frac{rq^{5r+2}}{1-q^{5r+2}} \right)+ \sum_{r \geq 0} \left(\frac{2(r+1)q^{5r+4}}{1-q^{5r+4}} - \frac{(r+1)q^{5r+3}}{1-q^{5r+3}}\right) \\- 2\sum_{r \geq 1} \frac{q^{5r}}{1-q^{5r}}. \tag{3.4}
\end{multline} 
\end{corollary}

Again, in \textbf{Proposition 2}, let $q \to q^5, b=c=q^2,$ to get :
\vspace{0.69cm}
\begin{corollary}
  \begin{multline} \label{eqMod5(b)}
    \left\{ \sum_{r \geq 0} \left( \frac{q^{5r+1}}{1-q^{5r+1}}+\frac{q^{5r+2}}{1-q^{5r+2}}-\frac{q^{5r+3}}{1-q^{5r+3}} -\frac{q^{5r+4}}{1-q^{5r+4}} \right) \right\}^2 \\  = \sum_{r \geq 0} \frac{q^{5r+1}}{(1-q^{5r+1})^2} +  \sum_{r \geq 0} \frac{q^{5r+4}}{(1-q^{5r+4})^2}+ \sum_{r \geq 1} \left( \frac{2rq^{5r+2}}{1-q^{5r+2}} - \frac{rq^{5r+4}}{1-q^{5r+4}} \right)+ \sum_{r \geq 0} \left(\frac{2(r+1)q^{5r+3}}{1-q^{5r+3}} - \frac{(r+1)q^{5r+1}}{1-q^{5r+1}}\right) \\- 2\sum_{r \geq 1} \frac{q^{5r}}{1-q^{5r}}. \tag{3.5}
\end{multline} 
\end{corollary}
Now, if we subtract \eqref{eqMod5(b)} from \eqref{eqMod5(a)}, we arrive at : 
\vspace{0.69cm}
\begin{corollary}
  \begin{multline} \label{eqMod5(c)}
    \left\{ \sum_{r \geq 0} \left( \frac{q^{5r+1}}{1-q^{5r+1}} -\frac{q^{5r+4}}{1-q^{5r+4}} \right) \right\}\left\{ \sum_{r \geq 0} \left( \frac{q^{5r+2}}{1-q^{5r+2}} -\frac{q^{5r+3}}{1-q^{5r+3}} \right) \right\} \\  = \frac{1}{4}\sum_{r \geq 0} \left( \frac{q^{5r+1}}{(1-q^{5r+1})^2} - \frac{q^{5r+2}}{(1-q^{5r+2})^2} - \frac{q^{5r+3}}{(1-q^{5r+3})^2} + \frac{q^{5r+4}}{(1-q^{5r+4})^2}\right)+ \frac{3}{4}\sum_{r \geq 1} \left( \frac{rq^{5r+2}}{1-q^{5r+2}} - \frac{rq^{5r+1}}{1-q^{5r+1}} \right) \\+ \frac{3}{4}\sum_{r \geq 1}\left( \frac{rq^{5r+3}}{1-q^{5r+3}} - \frac{rq^{5r+4}}{1-q^{5r+4}} \right)  -\frac{1}{4} \sum_{r \geq 0} \frac{q^{5r+1}}{1-q^{5r+1}} -\frac{1}{2} \sum_{r \geq 0} \frac{q^{5r+4}}{1-q^{5r+4}} + \frac{3}{4} \sum_{r \geq 0} \frac{q^{5r+3}}{1-q^{5r+3}}. \tag{3.6}
\end{multline} 
\end{corollary}

Now we move on to examine \textbf{Proposition 3}. Let $ q \to q^7$, $b=q$, $c=q^2$, then we have : 
\vspace{0.69cm}
\begin{corollary}
  \begin{flalign*}
    \left\{ \frac{1}{2} + \sum_{r \geq 0} \left( \frac{q^{7r+1}}{1-q^{7r+1}}+\frac{q^{7r+2}}{1-q^{7r+2}}-\frac{q^{7r+3}}{1-q^{7r+3}} +\frac{q^{7r+4}}{1-q^{7r+4}} -\frac{q^{7r+5}}{1-q^{7r+5}} -\frac{q^{7r+6}}{1-q^{7r+6}}\right) \right\}^2 \\ \tag{3.7} =\frac{1}{4} + \sum_{r \geq 1} \frac{rq^r}{1-q^r} - \sum_{r \geq 1} \frac{7rq^{7r}}{1-q^{7r}}.
\end{flalign*}  
\end{corollary}
As our last insight, let us consider the following expression :
\begin{multline} 
\left\{\frac{1}{2}+l(a)+l(b)-l(ab)\right\}^2 - \left\{\frac{1}{2}+l(c)+l(d)-l(cd)\right\}^2 \\ \notag = \left\{l(a)+l(b)-l(c)-l(d)-l(ab)+l(cd)\right\}\left\{ 1+l(a)+l(b)+l(c)+l(d)-l(ab)-l(cd) \right\}.
\end{multline}

Since the expression from the left hand side can be resolved in view of \textbf{Proposition 3}, we have : 
\begin{multline}
    \left\{l(a)+l(b)-l(c)-l(d)-l(ab)+l(cd)\right\}\left\{ 1+l(a)+l(b)+l(c)+l(d)-l(ab)-l(cd) \right\} \\=  \sum_{r \geq 0} \left\{ \frac{aq^r}{(1-aq^r)^2} +  \frac{bq^r}{(1-bq^r)^2} +\frac{abq^r}{(1-abq^r)^2}-\frac{cq^r}{(1-cq^r)^2} -  \frac{dq^r}{(1-dq^r)^2} -\frac{cdq^r}{(1-cdq^r)^2} \right\} \\+ \sum_{r \geq 1} \left\{ \frac{a^{-1}q^r}{(1-a^{-1}q^r)^2} +  \frac{b^{-1}q^r}{(1-b^{-1}q^r)^2} +\frac{a^{-1}b^{-1}q^r}{(1-a^{-1}b^{-1}q^r)^2} -\frac{c^{-1}q^r}{(1-c^{-1}q^r)^2} -  \frac{d^{-1}q^r}{(1-d^{-1}q^r)^2} -\frac{c^{-1}d^{-1}q^r}{(1-c^{-1}d^{-1}q^r)^2}\right\}  \tag{3.8} \label{12TermEq} 
\end{multline} 
Now we end our exposition with the following corollary. Put $q \to q^{13}$ and $a=q, b=q^3, c=q^2, d=q^6,$ in \eqref{12TermEq} to arrive at :
\vspace{0.69cm}
\begin{corollary} For $\chi_1, \chi_2, \chi_3$ defined as :

\begin{flalign*}
    \chi_1(n) = \begin{cases}
        0 & \text{\textit{if }} n \equiv 0 \pmod{13},\\
        1 & \text{\textit{if }} n \equiv 1, 3, 7, 8, 9, 11 \pmod{13},\\
        -1 & \text{\textit{if }} n \equiv 2, 4, 5, 6, 10, 12 \pmod{13}.
    \end{cases}\\
    \chi_2(n) = \begin{cases}
        0 & \text{\textit{if }} n \equiv 0 \pmod{13},\\
        1 & \text{\textit{if }} n \equiv 1, 2, 3, 5, 6, 9 \pmod{13},\\
        -1 & \text{\textit{if }} n \equiv 4, 7, 8, 10, 11, 12 \pmod{13}.
    \end{cases}\\
    \chi_3(n) = \begin{cases}
        0 & \text{\textit{if }} n \equiv 0 \pmod{13},\\
        1 & \text{\textit{if }} n \equiv 1, 3, 4, 9, 10, 12 \pmod{13},\\
        -1 & \text{\textit{if }} n \equiv 2, 5, 6, 7, 8, 11 \pmod{13}.
        \end{cases}
\end{flalign*}
There holds : 
  \begin{flalign*}
        \tag{3.9}\left\{ \sum_{r \geq 1}  \frac{\chi_1(r)q^r}{1-q^r} \right\}\left\{ 1+\sum_{r \geq 1}  \frac{\chi_2(r)q^r}{1-q^r}   \right\} = \sum_{r \geq 1}  \frac{\chi_3(r)q^r}{(1-q^r)^2} .
\end{flalign*}  
\end{corollary}

\printbibliography
\end{document}